\documentclass[twoside,11pt]{amsart}

\usepackage{amsmath,latexsym,amssymb,mathptm, times}
\input amssym.def
\input amssym
\input xypic
\input xy
\xyoption{all}
\def\demo{\noindent{\bf Proof. }}
\def\sqr#1#2{{\vcenter{\hrule height.#2pt
        \hbox{\vrule width.#2pt height#1pt \kern#1pt
                \vrule width.#2pt}
        \hrule height.#2pt}}}
\def\square{\mathchoice\sqr64\sqr64\sqr{4}3\sqr{3}3}
\def\QED{\hfill$\square$}

\def\m{{\mathfrak m}}
\def\n{{\mathfrak n}}

\def\p{{\mathfrak p}}
\def\ms{\medskip}

\newtheorem{Theorem}{Theorem}[section]
\newtheorem{Lemma}[Theorem]{Lemma}

\newtheorem{Proposition}[Theorem]{Proposition}

\newtheorem{Notation and Discussion}[Theorem]{Notation and Discussion}
\newtheorem{Assumptions and Discussion}[Theorem]{Assumptions and Discussion}

\newtheorem{Remark}[Theorem]{Remark}
\newtheorem{Example}[Theorem]{Example}

\begin{document}

\baselineskip=16pt

\title[Generalized Hilbert functions]
{\Large\bf Generalized Hilbert functions}

\author[C. Polini and Y. Xie]
{Claudia Polini \and Yu Xie}

\thanks{AMS 2010 {\em Mathematics Subject Classification}.
Primary 13D40 and 13A30; Secondary 13H15.}

\thanks{The first author was partially supported by the
NSF and the NSA}

\address{Department of Mathematics, University of Notre Dame,
Notre Dame, Indiana 46556} \email{cpolini@nd.edu}

\address{Department of Mathematics,University of Notre Dame,
Notre Dame, Indiana 46556} \email{yxie@nd.edu}

\vspace{-0.1in}

\begin{abstract}
Let $M$ be a finite module and let $I$ be an arbitrary ideal over a Noetherian local ring. We define the  generalized Hilbert  function of $I$ on $M$
using the $0$th local cohomology functor.  We show that our definition re-conciliates with that of Ciuperc$\breve{{\rm a}}$. By  generalizing Singh's formula (which holds in the case of $\lambda(M/IM)<\infty$), we prove that the generalized Hilbert coefficients $j_0, \ldots, j_{d-2}$ are preserved under a general hyperplane section, where $d={\rm dim}\,M$. We also keep track of the  behavior of $j_{d-1}$. Then we apply these results to study  the generalized Hilbert function for ideals that have minimal $j$-multiplicity or  almost minimal  $j$-multiplicity.  We provide counterexamples to show that  the  generalized Hilbert series of ideals having minimal or almost minimal  $j$-multiplicity does not have the `expected' shape described in the case where $\lambda(M/IM)<\infty$. Finally we  give a sufficient  condition such that the
generalized Hilbert series has the desired  shape.
\end{abstract}

\maketitle

\vspace{-0.2in}

\section{Introduction}

\bigskip

In this paper we define  the generalized Hilbert function   of  $I$, where $I$ is an arbitrary ideal over a Noetherian local ring $(R, \m)$. We study its  behavior under general hyperplane sections and  investigate the interplay between the generalized Hilbert function and the depth of the associated graded ring. The associated graded ring ${\rm gr}_I(R):=\oplus_{t=1}^{\infty}I^t/I^{t+1}$ of $I$ is an
algebraic construction whose projective scheme represents the
exceptional fiber of the blowup of a variety along a subvariety. The Hilbert function of $I$ provides useful information on its
arithmetical properties, like its depth, which in turn  give  information,
for instance, on the cohomology groups of the blowup.

In this introduction we will only discuss the case of associated graded rings,
although in the rest of the paper we will treat associated graded modules.

For an $\m$-primary ideal $I$, the Hilbert-Samuel function  is defined to be the numerical function $H_I(t)$ that measures the growth of the length $\lambda(R/I^{t+1})$ of the $(t+1)$th power of $I$ for all $t\geq 0$. For $t$ sufficiently large, the function $H_I(t)$ is a polynomial in $t$ of degree $d$, the dimension of $R$. This is the Hilbert-Samuel polynomial $P_I(t)$ of $I$, whose normalized coefficients $e_i(I)$, $0\leq i\leq d$, dubbed the Hilbert-Samuel coefficients of $I$, are uniquely determined by $I$. The notion  of Hilbert function is  fundamental and  has been widely studied  by algebraists and geometers. For example, the well-known Singh's formula  proved for the Hilbert-Samuel  function yields the fact
that the   Hilbert-Samuel coefficients $e_0, \ldots, e_{d-2}$ are preserved  under a general  hyperplane section. It also allows us to
keep track of  the  behaviors of $e_{d-1}$ and $e_{d}$ under a general  hyperplane section.


The interplay between the Hilbert function
of $I$, more precisely its Hilbert coefficients, and the depth
of the  associated graded ring  has been widely investigated. This
line of study has its roots in the pioneering work of Sally.
The idea is that extremal values of the Hilbert
coefficients, most
notably of the multiplicity of $I$,  yield high depth of the
associated graded ring and, conversely, good depth properties encode
all the information about the Hilbert function.
In 1967 Abhyankar proved that the multiplicity  $e(R)=e_0(\m)$ of a
$d$-dimensional Cohen-Macaulay local ring is bounded below by
$\mu(\m)-d+1$, where $\mu(\m)$ is the embedding dimension of $R$
\cite{AB}. Rings for which $e_0=\mu(\m)-d+1$ have since then been
called rings of minimal multiplicity. In the case of minimal
multiplicity, Sally had shown   that the associated
graded ring ${\rm gr}_\m(R)$ is always Cohen-Macaulay with the Hilbert-Samuel
 series $h_I(z):=\sum_{t=0}^{\infty}H_I(t)z^t=(1+hz)/(1-z)^{d+1}$, 
 where $h=\mu(\m)-d$ is the embedding codimension of $R$ \cite{S1}.
Almost twenty years later, Rossi-Valla \cite{RV1} and   Wang \cite{W} independently proved that if the
multiplicity of $R$ is almost minimal then the depth of ${\rm gr}_\m(R)$ is almost maximal, i.e., it is at least
$d-1$. 
Furthermore in \cite{RV1}, all the possible Hilbert-Samuel series of rings 
of almost minimal multiplicity have been described to 
be of the form $h_I(z)=(1+hz+z^s)/(1-z)^{d+1}$, where $2\leq s\leq h+1$.
Since then  there have been many generalizations of
these results to $\m$-primary ideals and to modules having ideal filtrations of finite colength,
 a condition that
is required to define the classical Hilbert function (see for
example \cite{RV2}, \cite{H}, \cite{CPV}, \cite{E}, \cite{R1},  \cite{P1}, \cite{RV}).

The first work to generalize the theory of minimal and almost minimal 
multiplicity to arbitrary ideals was done
 by Polini and Xie \cite{PX}.
If the ideal $I$ is not $\m$-primary then the Hilbert function is not defined,
thus there is no numerical information on Hilbert coefficients
available to study the Cohen-Macaulayness of ${\rm gr}_I(R)$. To
remedy the lack of this tool, Polini and Xie  proposed to use the
notion of $j$-multiplicity. The $j$-multiplicity was  introduced by Achilles
and Manaresi \cite{AM} in 1993 as a generalization of the Hilbert multiplicity to arbitrary
ideals.
It has been frequently used by both
algebraists and geometers as an invariant to deal with improper
intersections and  non-isolated singularities
\cite{AM}.

In \cite{PX}, Polini and Xie  proved a lower bound for the $j$-multiplicity of
an ideal  and then gave  the definition of ideals having minimal
$j$-multiplicity or almost minimal $j$-multiplicity. Under certain residual
conditions, they proved that for  ideals having minimal
$j$-multiplicity (respectively, almost minimal $j$-multiplicity), the associated graded ring
is Cohen-Macaulay (respectively, almost Cohen-Macaulay).
Their work has been extended  further very recently by Mantero and Xie by introducing the concept of $j$-stretched ideals \cite{MX}.
Although  both work successfully generalized many classical results,   there is no description of the Hilbert function as it has been done for  $\m$-primary ideals.

In this paper, we define the generalized Hilbert function of an arbitrary ideal $I$ on a finite module $M$ over a Noetherian local ring  using 0th local cohomology. We then present  Ciuperc$\breve{{\rm a}}$'s approach and re-conciliate that with ours. We generalize Singh's formula to arbitrary ideals and show that the first $d-1$ generalized Hilbert coefficients $j_0, \ldots, j_{d-2}$ are preserved under a general hyperplane section, where $d$ is the dimension of the module $M$. Furthermore, we keep track of how the next generalized Hilbert coefficient 
$j_{d-1}$ changes. In Section 3, we study
the generalized Hilbert function for ideals that have minimal $j$-multiplicity or have almost minimal  $j$-multiplicity   on $M$.
 We provide counterexamples to show that the generalized Hilbert series of ideals having minimal or almost minimal   $j$-multiplicity on $M$ does not have the `expected' shape as described in the case where $\lambda(M/IM)<\infty$. We then give a sufficient  condition such that the generalized Hilbert series has the desired shape.


The study of the generalized Hilbert function  is of independent interest and for example will be used in a future work where we investigate the first generalized Hilbert  coefficient $j_1$ in relation with the normalization of the Rees algebra
\cite{NPUX}.

\section{Basic properties of generalized Hilbert functions}

\medskip


Let $M$ be a finite module over a Noetherian local ring $(R, \m)$ and  $I$  an arbitrary
$R$-ideal. Let $T={\rm gr}_{I} (M):=\oplus_{t=0}^{\infty}
I^tM/I^{t+1}M$ be the {\it associated graded module}  of $I$ on $M$. Observe $T$ is a finite graded module over the {\it associated graded ring}
 $G={\rm
gr}_I(R):=\oplus_{t=0}^{\infty}I^t/I^{t+1}$.  In general the homogeneous components of  $T$ may
not have  finite length, thus we consider the $T$-submodule of
elements supported on $\m$,  $W=\Gamma_{\m}(T)=0:_T
{\m}^{\infty}=\oplus_{t=0}^{\infty} \Gamma_{\m}\,(I^tM/I^{t+1}M)$.
Since $W$ is annihilated by a large power of $\m$, it is a finite
graded module over ${\rm gr}_I(R) \otimes_R R/{\m}^{\alpha}$ for some $\alpha\geq 0$,
hence its Hilbert-Samuel  function $H_W(t)=\sum_{\nu=0}^t \lambda(\Gamma_{\m}\,(I^\nu M/I^{\nu +1}M))$ is well defined. We define the {\it generalized Hilbert-Samuel function} of $I$ on $M$: $H_{I,\, M}(t):=H_W(t)$ for every $t\geq 0$, and the {\it generalized Hilbert-Samuel series} of $I$ on $M$: $h_{I,\, M}(z):=\sum_{t=0}^{\infty}H_W(t)z^t$. By the Hilbert-Serre theorem, we know that the series is of the form
$$
h_{I,\, M}(z)=h(z)/(1-z)^{r+1},
$$
where $h(z)\in \mathbb{Z}[z]$, $h(1)\neq 0$ and $r={\rm dim}_G\,W$. The polynomial $h(z)$ is called the {\it $h$-polynomial} of $I$ on $M$.


Let ${\rm dim}\,M=d$. Then ${\rm dim}_G\,W\leq {\rm dim}_G\,T=d$, thus  $H_{I,\, M}(t)$ is eventually a polynomial $$P_{I,\, M}(t)=\sum_{i=0}^{d} (-1)^i j_i(I, M)
{{t+d-i}\choose{d-i}}$$ of degree  at most $d$. We define $P_{I,\, M}(t)$ to be the {\it generalized Hilbert-Samuel polynomial} of $I$ on $M$ and $j_i(I, M),\, 0\leq i\leq d$, the {\it generalized Hilbert-Samuel coefficients} of $I$ on $M$. Observe that  the normalized leading coefficient $j_0(I, M)$ is also called the {\it $j$-multiplicity} of $I$   on  $M$ (see for instance \cite{AM} or \cite{NU}).

Recall that the Krull dimension of the special fiber module
$T/{\m}T$ is called the
 {\it analytic spread} of $I$ on $M$ and is denoted by $\ell(I, M)$. In general,
 ${\rm dim}_G\,W \leq \ell(I,M) \leq d$
 and equalities hold if and only if $\ell(I,M)=d$. Therefore $j_0(I, M)\neq 0$ if and only if $\ell(I, M)=d$ \cite[2.1]{NU}.

If $M/IM$ has finite length, the ideal $I$ is said to be  an {\it
ideal of definition} on $M$. In this case each homogeneous
component of $T$ has finite length, thus $W=T$ and the
generalized Hilbert-Samuel function coincides with the usual Hilbert-Samuel function; in particular the generalized Hilbert-Samuel coefficients $j_i(I,M)$ coincides with the usual Hilbert-Samuel coefficients $e_i(I,M)$, $0\leq i\leq d$.



The definition of generalized Hilbert function, generalized Hilbert polynomial, and corresponding generalized Hilbert coefficients is new to the literature. It appeared first in an unpublished work of Polini and Ulrich. It stems from the original definition of $j$-multiplicity given by Achilles and Manaresi \cite{AM}. The definition of  higher generalized Hilbert coefficients presented in the literature is quite different (see \cite{AM1} and \cite{C}). We review now Ciuperc$\breve{{\rm a}}$'s approach and re-conciliate that with ours. More precisely, for each $i$, Ciuperc$\breve{{\rm a}}$ defined a sequence of numbers which represents the $i$th generalized Hilbert coefficient \cite[2.1]{C}. We will show that our $(-1)^i j_i(I,M)$ is the last element in that sequence (see Proposition \ref{j-coefficientsBigraded}).

Let $q$ be an ideal of definition on $M/IM$. We form  the bigraded module $T^{\prime}={\rm gr}_q({\rm gr}_I(M))=\oplus_{s,t=0}^{\infty}(q^sI^tM+I^{t+1}M)/(q^{s+1}I^tM+I^{t+1}M)$, then  $T^{\prime}$ is a finite bigraded module over the bigraded ring $G^{\prime}=\oplus_{s,t=0}^{\infty}(q^sI^t+I^{t+1})/(q^{s+1}I^t+I^{t+1})$. Notice that every homogeneous component of $T^{\prime}$ has finite length. Let
$$
H_{(q,\, I,\, M)}^{(0,\,0)}(s,t)=
\lambda((q^sI^tM+I^{t+1}M)/(q^{s+1}I^tM+I^{t+1}M)),
$$
$$
H_{(q,\, I,\, M)}^{(1,\,0)}(s,t)=\sum_{\mu=0}^{s}H_{(q,\, I,\, M)}^{(0,\,0)}(\mu,t)
=\lambda(I^tM/(q^{s+1}I^tM+I^{t+1}M)),
$$
and
$$
H_{(q,\, I,\, M)}^{(1,\,1)}(s,t)=\sum_{\nu=0}^{t}H_{(q,\, I,\, M)}^{(1,\,0)}(s,\nu)=
\lambda(M/(q^{s+1}M+IM))
+\cdots+\lambda(I^tM/(q^{s+1}I^tM+I^{t+1}M)).
$$

\medskip

For  $s$ and $ t$  sufficiently large, $H_{(q,\, I,\, M)}^{(0,\,0)}(s,t)$, $H_{(q,\, I,\, M)}^{(1,\,0)}(s,t)$ and $H_{(q,\, I,\, M)}^{(1,\,1)}(s,t)$ are  polynomials $P_{(q,\, I,\, M)}^{(0,\,0)}(s,t)$, $P_{(q,\, I,\, M)}^{(1,\,0)}(s,t)$ and $P_{(q,\, I,\, M)}^{(1,\,1)}(s,t)$ in two variables $s, t$ of degree at most $d-2$, $d-1$ and $d$ respectively. One can write
$$P_{(q,\, I,\, M)}^{(0,\,0)}(s,t)=\sum_{i,\,j\geq 0,\, i+j\leq d-2}a_{(q,\, I,\, M)}^{(0,\,0)}(i,j)\bigg(\begin{array}{c}s+i\\i \end{array}\bigg)\bigg(\begin{array}{c}t+j\\j \end{array}\bigg),
$$
$$P_{(q,\, I,\, M)}^{(1,\,0)}(s,t)=\sum_{i,\,j\geq 0,\, i+j\leq d-1}a_{(q,\, I,\, M)}^{(1,\,0)}(i,j) \bigg(\begin{array}{c}s+i\\i \end{array}\bigg)\bigg(\begin{array}{c}t+j\\j \end{array}\bigg),
$$
and
$$
P_{(q,\, I,\, M)}^{(1,\,1)}(s,t)=\sum_{i,\,j\geq 0,\, i+j\leq d}a_{(q,\, I,\, M)}^{(1,\,1)}(i,j)\bigg(\begin{array}{c}s+i\\i \end{array}\bigg)\bigg(\begin{array}{c}t+j\\j \end{array}\bigg).
$$
\medskip
Then $a_{(q,\, I,\, M)}^{(0,\,0)}(i,j)=a_{(q,\, I,\, M)}^{(1,\,0)}(i+1,j)$ for $i,\,j\geq 0,\, i+j\leq d-2$ and $a_{(q,\, I,\, M)}^{(1,\,0)}(i,j)=a_{(q,\, I,\, M)}^{(1,\,1)}(i,j+1)$ for $i,\,j\geq 0,\, i+j\leq d-1$.
The $i$th {\it generalized Hilbert coefficient} defined by  Ciuperc$\breve{{\rm a}}$   is  as follows \cite[2.1]{C}:
$$
j_i(q, I, M)=(a_{(q,\, I,\, M)}^{(1,\,1)}(i,d-i),\, a_{(q,\, I,\, M)}^{(1,\,1)}(i-1,d-i), \ldots, a_{(q,\, I,\, M)}^{(1,\,1)}(0,d-i))\in \mathbb{Z}^{i+1}
$$
for every $0\leq i \leq d$.

Notice that  $j_0(q, I, M)=a_{(q,\, I,\, M)}^{(1,\,1)}(0,d)=a_{(q,\, I,\, M)}^{(1,\,0)}(0,d-1)$ and $a_{(q,\, I,\, M)}^{(1,\,1)}(0,d)=j_0(I, M)$  the $j$-multiplicity of $I$ on $M$ (see \cite{AM1} or \cite{C}). Moreover, if  $I$ is an ideal of definition on $M$, then for $0\leq i \leq d$, $a_{(q,\, I,\, M)}^{(1,\,1)}(0,d-i)=(-1)^i j_i(I, M)=(-1)^i e_i(I,M)$  and $a_{(q,\, I,\, M)}^{(1,\,1)}(l,d-i)=0$ if $0<l\leq i$ (\cite{AM1}, \cite{C}). In particular,
$a_{(q,\, I,\, M)}^{(1,\,1)}(0,d)=j_0(I, M)=e_0(I, M)$ is the usual Hilbert-Samuel multiplicity of $I$ on $M$. Our first goal is to show that for any ideal $I$, the coefficient $a_{(q,\, I,\, M)}^{(1,\,1)}(0,d-i)$ coincides  with $(-1)^i j_i(I, M)$ for suitable choice of $q$ (see Proposition \ref{j-coefficientsBigraded}).
\ms

In the following lemma (Lemma \ref{j-multiplicityformula}) we are going to show that Ciuperc$\breve{{\rm a}}$'s
$i$th generalized Hilbert coefficient behaves well under a general hyperplane section for $0\leq i \leq d-2$.
This generalizes  Proposition 2.11 in \cite{C}
since we do not assume that the hyperplane section is a nonzero divisor.
Furthermore, we keep track of how the $(d-1)$th generalized Hilbert coefficient  
$j_{d-1}(q, I, M)$ changes (notice that there is an error in \cite[2.11]{C}, namely $j_{d-1}(q, I, M)$ is not preserved, see Example \ref{GeneralizedCoefficient}).
The proof uses a generalization of Singh's formula (see Lemma \ref{Singh})
as it is done in the case of ideals of definition on $M$  (see \cite[Lemma 1.6, Proposition 1.2]{RV}).


\begin{Lemma}\label{Singh} Let $R$ be a Noetherian local ring, let $I$ be any $R$-ideal, and let $M$ be a finitely generated $R$-module.  For any $x\in I$ write $\overline{M}=M/xM$. The generalized Singh's formula holds for any non negative integers $s$ and $t$:
$$\vspace{0.07in}
\begin{array}{l}
H_{(q,\, I,\, M)}^{(1,\,0)}(s,t)=H_{(q,\, I,\, \overline{M})}^{(1,\,1)}(s,t)\vspace{0.1in}\\
+\sum_{\nu=2}^t\lambda(I^{\nu}M:_M x/
[(q^{s+1}I^{\nu}M+I^{\nu+1}M):_M x +I^{\nu-1}M])\vspace{0.1in}\\
-\sum_{\nu=1}^t \lambda([(q^{s+1}I^\nu M+I^{\nu+1}M):_{I^{\nu-1}M}x]/[q^{s+1}I^{\nu-1}M+I^\nu M]).
\end{array}\vspace{0.07in}
$$
Moreover if $I$ is an ideal of definition on $M$, then for $s>>0$, the above formula becomes the usual Singh's formula:
\vspace{-0.1in}
$$
\lambda (I^{t}M/I^{t+1}M)=\lambda (M/I^{t+1}M+xM)-\lambda(I^{t+1}M:_M x/I^{t}M).
$$
\end{Lemma}

\demo
First, for every $\nu\geq 0$, there is an exact sequence:
$$ \vspace{0.07in}
\begin{array}{c}
0\longrightarrow [(q^{s+1}I^\nu M+I^{\nu+1}M):_{I^{\nu-1}M}x]/[q^{s+1}I^{\nu-1}M+I^\nu M]\longrightarrow \vspace{0.1in}\\ I^{\nu-1}M/[q^{s+1}I^{\nu-1}M+I^\nu M] \stackrel{x}{\longrightarrow} I^{\nu}M/[q^{s+1}I^{\nu}M+I^{\nu+1}M] \vspace{0.1in}\\ \longrightarrow I^{\nu}M/[xI^{\nu-1}M+q^{s+1}I^{\nu}M+I^{\nu+1}M]\longrightarrow 0.
\end{array}\vspace{0.1in}
$$
By adding $\nu$ up to $t$,  we obtain that
\begin{align}
&H_{(q,\, I,\, M)}^{(1,\,0)}(s,t)=H_{(q,\, I,\, \overline{M})}^{(1,\,1)}(s,t) \notag\\
&+\sum_{\nu=2}^t \lambda ([xM\cap I^{\nu}M+q^{s+1}I^{\nu}M+I^{\nu+1}M]/[xI^{\nu-1}M+q^{s+1}I^{\nu}M+I^{\nu+1}M]) \notag\\
&-\sum_{\nu=1}^t \lambda([(q^{s+1}I^\nu M+I^{\nu+1}M):_{I^{\nu-1}M}x]/[q^{s+1}I^{\nu-1}M+I^\nu M]).
\end{align}
Observe that
$$
\begin{array}{l}
[xM\cap I^{\nu}M+q^{s+1}I^{\nu}M+I^{\nu+1}M]/[xI^{\nu-1}M+q^{s+1}I^{\nu}M+I^{\nu+1}M] \vspace{0.1in}\\
\cong  xM\cap I^{\nu}M/[xM \cap (q^{s+1}I^{\nu}M+I^{\nu+1}M)+xI^{\nu-1}M] \vspace{0.1in}\\
= x(I^{\nu}M:_M x)/
x[(q^{s+1}I^{\nu}M+I^{\nu+1}M):_M x +I^{\nu-1}M].
\end{array}
$$
It is easy to see that the surjective map $I^{\nu}M:_M x \stackrel{x}{\rightarrow} x(I^{\nu}M:_M x)/
x[(q^{s+1}I^{\nu}M+I^{\nu+1}M):_M x +I^{\nu-1}M]$ has kernel $(q^{s+1}I^{\nu}M+I^{\nu+1}M):_M x +I^{\nu-1}M$. Therefore
$$ \vspace{0.1in}
\begin{array}{l}
\lambda ([xM\cap I^{\nu}M+q^{s+1}I^{\nu}M+I^{\nu+1}M]/[xI^{\nu-1}M+q^{s+1}I^{\nu}M+I^{\nu+1}M])\vspace{0.1in}\\
=\lambda (x(I^{\nu}M:_M x)/
x[(q^{s+1}I^{\nu}M+I^{\nu+1}M):_M x +I^{\nu-1}M])\vspace{0.1in}\\
=\lambda(I^{\nu}M:_M x/
[(q^{s+1}I^{\nu}M+I^{\nu+1}M):_M x +I^{\nu-1}M]).
\end{array}
$$
By plugging in equation (1) we obtain the desired formula.

For the second part, observe that if $I$ is an ideal of definition on $M$ then for $s>>0$, $q^{s+1}M\subseteq IM$. Hence  the  formula in the statement of Lemma \ref{Singh}  becomes
\begin{eqnarray*}
& &\lambda (I^{t}M/I^{t+1}M)-\lambda (M/I^{t+1}M+xM) \\
&=& \sum_{\nu=2}^t \lambda(I^{\nu}M:_M x/
[I^{\nu+1}M:_M x +I^{\nu-1}M]) -\sum_{\nu=1}^t \lambda(I^{\nu+1}M:_{I^{\nu-1}M}x/I^\nu M)
\\ &=& \sum_{\nu=2}^t \lambda(I^{\nu}M:_M x/
I^{\nu-1}M)- \sum_{\nu=2}^t \lambda([I^{\nu+1}M:_M x+I^{\nu-1}M]/
I^{\nu-1}M)-\sum_{\nu=1}^t \lambda(I^{\nu+1}M:_{I^{\nu-1}M}x/I^\nu M)
\\ &=& \sum_{\nu=2}^t \lambda(I^{\nu}M:_M x/
I^{\nu-1}M)- \sum_{\nu=2}^t \lambda(I^{\nu+1}M:_M x/
I^{\nu+1}M:_{I^{\nu-1}M} x)-\sum_{\nu=1}^t \lambda(I^{\nu+1}M:_{I^{\nu-1}M}x/I^\nu M)
\\ &=& \sum_{\nu=2}^t \lambda(I^{\nu}M:_M x/
I^{\nu-1}M)- \sum_{\nu=1}^t \lambda(I^{\nu+1}M:_M x/
I^{\nu}M)\\
&=&-\lambda(I^{t+1}M:_M x/I^tM).
\end{eqnarray*}
\QED

\smallskip

We now need to recall some  definitions and basic facts (see for instance \cite{AM}, \cite{Hu}, \cite{Trung} and \cite{PX}).
Let $(R, \m)$ be a Noetherian local ring, $I$ an $R$-ideal, $M$ a finite $R$-module, and $G$ and $T$ are defined as before.
An element $x\in R$ is said to be a {filter-regular element} for $M$ with respect to $I$, if $0:_M x\subseteq 0:_M I^{\infty}$. 
This holds if and only if $x$ avoids all associated primes of $M$ that do not contain $I$.
 A sequence of elements $x_1, \ldots, x_{l}$ of $I$ is called a {\it filter-regular sequence} for $M$ with respect to $I$, if $(x_1, \ldots, x_{i-1})M:_M x_i\subseteq (x_1, \ldots, x_{i-1})M:_M I^{\infty}$ for $1\leq i\leq l$. 
 The sequence  $x_1, \ldots, x_{l}$ of $I$ is called a {\it superficial sequence} for $M$ with respect to $I$ if the initial forms $x_1^*,\ldots, x_{l}^*$ of $x_1, \ldots, x_{l}$ in $G$ are of degree one and form a filter-regular sequence for $T$ with respect to $G_+$, where $G_+$ is the ideal generated by all homogeneous elements of positive degree in $G$. It is well-known that if $x_1, \ldots, x_{l}$ form  a {\it superficial sequence} then  they also form a  filter-regular sequence.
 
Let $I=(a_1, \ldots, a_n)$.  Recall that $a_1, \ldots, a_n$ form a {\it $d$-sequence} on $M$ (or {\it an absolutely superficial $M$-sequence} in the sense of Trung) if $[(a_1, \ldots, a_{i-1})M:_M a_{i}]\cap IM=(a_1, \ldots, a_{i-1})M$ for $1\leq i\leq n$.

Let $I=(a_1, \ldots, a_n)$ and write $x_i=\sum_{j=1}^n
\lambda_{ij}a_j$ for $1\leq i\leq l$ and $(\lambda_{ij})\in R^{l
n}$.  The elements  $x_1, \ldots, x_{l}$ form a {\it sequence of
general elements} in $I$ (equivalently  $x_1, \ldots, x_{l}$ are {\it
general} in $I$) if there exists a dense open subset $U$ of $K^{l n}$
such that the image $(\overline{\lambda_{ij}})\in U$, where $K=R/\m$ is the residual field of $R$.
When $l=1$,  $x=x_1$ is said to be {\it general} in $I$.

The notion of general elements is a fundamental tool for our study
as they are always a superficial sequence for $M$ with respect to $I$ \cite[2.5]{X}; 
they generate a minimal reduction $J$ whose reduction number
$r_J(I, M)$ coincides with the reduction number $r(I,M)$ of $I$ on
$M$ if $l=\ell(I, M)$ (see \cite[2.2]{T2} and \cite[8.6.6]{SH}); and they form a {\it super-reduction} in the sense of \cite{AM} whenever
$l=\ell(I, M)=d={\rm dim}_R\,M$ (see \cite[2.5]{X}).  Furthermore, one can compute the $j$-multiplicity using general
elements.

The next proposition states the behaviors of Ciuperc$\breve{{\rm a}}$'s generalized Hilbert coefficients under a general hyperplane section.

\begin{Proposition}\label{j-multiplicityformula}
Let $R$ be a Noetherian local ring with infinite residue field, let $I$ be any $R$-ideal, and let $M$ be a finitely generated $R$-module of dimension $d$.   If $x$ is a general element in $I$ then the following statements hold:
\begin{itemize}
\item[(a)] $j_i(q, I, M)=j_i(q, I, M/xM)$ for $0\leq i\leq d-2$.\vspace{0.04in}

\item[(b)] There exists a fixed positive integer $t_0$ such that if $t> t_0$,
$$\sum_{\nu=2}^{t}\lambda(I^{\nu}M:_M x/
[(q^{s+1}I^{\nu}M+I^{\nu+1}M):_M x +I^{\nu-1}M])
$$
$$
-\sum_{\nu=1}^{t} \lambda([(q^{s+1}I^\nu M+I^{\nu+1}M):_{I^{\nu-1}M}x]/[q^{s+1}I^{\nu-1}M+I^\nu M])
$$
is a polynomial $\sum_{i=0}^{d-1}\beta_i\bigg(\begin{array}{c}s+i\\i \end{array}\bigg)$ in $s$  of degree at most $d-1$ (i.e., it does not depend on $t$).
Therefore 
$j_{d-1}(q, I, M)=j_{d-1}(q, I,M/xM)+(\beta_{d-1}, \ldots, \beta_0).$
\end{itemize}
\end{Proposition}

\demo 
By avoiding finitely many associated  prime ideals of $T^{\prime}$ which do not contain $I/q I$, we obtain a general element $x$ in $I$ such that the image $x^{\prime}$ in $I/q I$  is filter-regular for $T^{\prime}$ with respect to the ideal   $(I/q I)G^{\prime}$.
Thus there exists $t_0 > 0$ such that for every $t> t_0$ and for every $s\geq 0$, \,
$xM\cap I^{t}M=xI^{t-1}M$ and $(q^{s+1}I^tM+I^{t+1}M):_{I^{t-1}M}x=q^{s+1}I^{t-1}M+I^tM$. Write $\overline{M}=M/xM$. By Lemma \ref{Singh} we obtain  for sufficiently large $s$ and $t$,
\begin{align}
P_{(q,\, I,\, M)}^{(1,\,0)}(s,t)=P_{(q,\, I, \, \overline{M})}^{(1,\,1)}(s,t) &+\sum_{\nu=2}^{t_0}\lambda(I^{\nu}M:_M x/
[(q^{s+1}I^{\nu}M+I^{\nu+1}M):_M x +I^{\nu-1}M]) \notag\\
&-\sum_{\nu=1}^{t_0} \lambda([(q^{s+1}I^\nu M+I^{\nu+1}M):_{I^{\nu-1}M}x]/[q^{s+1}I^{\nu-1}M+I^\nu M]).
\end{align}
Now the desired equality follows from  (2) by equaling the coefficients.
\QED

\vspace{0.3in}

The following lemma follows from \cite[1.1, 2.1]{Trung}.

\begin{Lemma}\label{multiplicity} Let $(R, \m)$ be a Noetherian local ring with infinite residue field, let $G$ be a standard graded ring over $R$, and let $T$ be  a finite  $G$-module.
Fix a positive integer $l$. Then there exists $\underline{y}=y_1, \ldots, y_{l}$ in $ \m$ such that $\underline{y}$\, forms a $d$-sequence on $T$, i.e.,  for every $1\leq i\leq l,$
$$(y_1, \ldots, y_{i-1})T:_T y_{i}=(y_1, \ldots, y_{i-1})T:_T \m^{\infty}
$$
and
$$
[(y_1, \ldots, y_{i-1})T:_T y_{i}]\cap (y_1, \ldots, y_{l})T=(y_1, \ldots,y_{i-1})T.
$$\vspace{0.05in}
In particular, if $T=\oplus_{t=0}^{\infty}T_t$  is a finite graded module over $G$, then for every $t\geq 0$,

$$(y_1, \ldots, y_{i-1})T_t:_{T_t} y_{i}=(y_1, \ldots, y_{i-1})T_t:_{T_t} \m^{\infty}$$ and
$$[(y_1, \ldots, y_{i-1})T_t:_{T_t} y_{i}]\cap (y_1, \ldots, y_{l})T_t=(y_1, \ldots,y_{i-1})T_t$$
which shows that $\underline{y}$\, forms an absolutely superficial $T_t$-sequence for every $t\geq 0$.
\end{Lemma}

\demo
Since $T$ is a finite $G$-module and $\m G$ is a $G$-ideal, there exist $a_1, \ldots, a_{l}$ in $\m$ such that their initial forms $a_1^*, \ldots, a_{l}^*$ in $\m G$ form a filter-regular sequence for $T$ with respect to $\m G$.
By \cite[2.1]{Trung}, there is an ascending sequence of integers $n_1\leq \n_2\leq \ldots \leq \n_{l}$ such that $(a_1^*)^{ n_1}, \ldots, (a_{l}^*)^{ n_{l}}$ form a $d$-sequence on $T$. Set $\underline{y}=a_1^{n_1}, \ldots, a_{l}^{n_{l}}$. Then the conclusions follow from \cite[1.1]{Trung}.
\QED

\vspace{0.2in}

In the following proposition we achieve our first goal. We show that how our generalized Hilbert coefficients are related to the generalized Hilbert coefficients of Ciuperc$\breve{{\rm a}}$. The proof of this result follows essentially from Lemma \ref{multiplicity} and \cite[4.1]{Trung}.


\begin{Proposition} \label{j-coefficientsBigraded}
Let $(R, \m)$ be a Noetherian local ring with infinite residue field, let $I$ be any $R$-ideal, let $M$ be a finitely generated $R$-module of dimension $d$, and  let $T$ be the associated graded module of $I$ on $M$. Then there exist elements $y_1, \ldots, y_d \in \m$ which form  a $d$-sequence on $T_t$ for every $t\geq 0$. Therefore if we set  $q=(y_1, \ldots, y_d)$ 
then 
$
a_{(q,\, I,\, M)}^{(1,\,1)}(0,d-i)=(-1)^{i}j_{i}(I, M)
$
for every $0\leq i \leq d$.
\end{Proposition}

\demo
By Lemma \ref{multiplicity}, there exist elements $y_1, \ldots, y_d \in \m$ such that $y_1, \ldots, y_d$ form   a $d$-sequence on $T_t$ for every $t\geq 0$, in particular, they form a system of parameters for $T_t$. By \cite[4.1]{Trung}, for every $t\geq 0$,
$$
H_{T^{\prime}}^{(1,0)}(s, t)=\lambda(T_t/q^{s+1}T_t)=\\
\sum_{i=0}^d [\lambda(q_{i}T_t:_{T_t} \m^{\infty}/q_{i}T_t)-\lambda(q_{i-1}T_t:_{T_t} \m^{\infty}/q_{i-1}T_t)]\left(\begin{matrix}
 s+i\\
 i\\
\end{matrix}\right ),
$$
where $q_i=(y_1, \ldots, y_i)$, $1\leq i\leq d$, $q_0=(0)$ and $q_{-1}T_t:_{T_t} \m^{\infty}/q_{-1}T_t=(0)T_t$.

Thus  for $t>>0$,
\begin{eqnarray*}
\sum_{i=0}^{d} (-1)^i j_i(I, M)
{{t+d-i}\choose{d-i}}&=&\sum_{\nu=0}^t\lambda(0:_{T_{\nu}} \m^{\infty})\\
&=& \sum_{i=0}^d a_{(q,\, I,\, M)}^{(1,\,1)}(0,d-i)
 {{t+d-i}\choose{d-i}}.
 \end{eqnarray*}
 We are done by equaling corresponding coefficients.
\QED

\vspace{0.2in}

We now recall the notion of specialization from \cite{T3}. Let $S=V[\underline{x}]$ be a polynomial ring over  a regular local ring $(V, \n)$   and  $\widetilde{S}=V[\underline{z}]_{\n V[\underline{z}]}[\underline{x}]$, where $\underline{z}=(z_1, \ldots, z_{l})$ are $l$ variables over $S$. Assume $V$ has infinite residue field $K$. Let $\underline{\alpha}=(\alpha_1, \ldots, \alpha_{l})$ be a family of elements of $V$. Let $\widetilde{F}$ be a free $\widetilde{S}$-module of finite rank. The specialization $\widetilde{F}_{\underline{\alpha}}$ is a free $S$-module  of the same  rank. Let $\widetilde{M}$ be a finite $\widetilde{S}$-module and $\phi: \widetilde{F}\rightarrow \widetilde{G}\rightarrow \widetilde{M}\rightarrow 0$  a finite free representation of $\widetilde{M}$ with representing matrix $B=(a_{ij}(\underline{z}, \underline{x}))$, where $a_{ij}(\underline{z}, \underline{x})=p_{ij}[\underline{z}, \underline{x}]/q[\underline{z}]$, $p_{ij}[\underline{z}, \underline{x}]\in V[\underline{z}, \underline{x}]$ and $q[\underline{z}]\in V[\underline{z}]\backslash \n V[\underline{z}]$. For any $\underline{\alpha}\in V^{l}$ such that $\overline{q[\underline{\alpha}]}\neq 0$ in $K$, we define $a_{ij}(\underline{\alpha}, \underline{x})=p_{ij}[\underline{\alpha}, \underline{x}]/q[\underline{\alpha}]$ and set $B_{\underline{\alpha}}=(a_{ij}(\underline{\alpha}, \underline{x}))$. Then $B_{\underline{\alpha}}$ is well defined for almost all $\underline{\alpha}$. The specialization $\phi_{\underline{\alpha}}: \widetilde{F}_{\underline{\alpha}}\rightarrow \widetilde{G}_{\underline{\alpha}}$ of $\phi$ is given by the matrix $B_{\underline{\alpha}}$ provided $B_{\underline{\alpha}}$ is well-defined. The specialization of $\widetilde{M}_{\underline{\alpha}}$ is defined as ${\rm coker}(\phi_{\underline{\alpha}})$. By \cite{T2}, $\widetilde{M}_{\underline{\alpha}}$ is well-defined and does not depend on the matrix representation $B$ for almost all $\underline{\alpha}$. It is easy to see that if $M$ is a finite $S$-module and  $\widetilde{M}=M\otimes_S \widetilde{S}$, then $\widetilde{M}_{\underline{\alpha}}\cong M$. Applying  a similar proof of \cite{T3}, one can show that all the properties of specialization proved in \cite{T3} when $V$ is a field also hold if $V$ is a regular local ring.

We arrive at our main theorem. It follows from  Propositions \ref{j-coefficientsBigraded} and \ref{j-multiplicityformula}.

\begin{Theorem}\label{j-coefficients}
Let $R$ be a Noetherian local ring with infinite residue field, let $I$ be any $R$-ideal, and let $M$ be a finitely generated $R$-module of dimension $d$.   If $x$ is a general element in $I$ then $ j_i(I,M)=j_i(I, M/xM)$ for $0\leq i\leq d-2$.
\end{Theorem}

\demo
By passing to the completion, we may assume that $R=V/H$, where $(V, \n)$ is a regular local ring.
We may further pass to $(V, \n)$ to assume that  $(R, \m)$ is a regular local ring. Write $I=(a_1, \ldots, a_n)$. Let $\widetilde{R}=R[\underline{z}]_{\m R[\underline{z}]}$, where $\underline{z}=(z_1, \ldots, z_n)$ are $n$ variables over $R$. Let $\widetilde{M}=M\otimes_R \widetilde{R}$, $\widetilde{T}={\rm gr}_{I\widetilde{R}}(\widetilde{M})=T\otimes_R \widetilde{R}$. Set $\widetilde{x}=\sum_{j=1}^nz_{j}a_j$ and $\widetilde{T}/\widetilde{x}\widetilde{M}:={\rm gr}_{I\widetilde{R}}(\widetilde{M}/\widetilde{x}\widetilde{M})$.
By Proposition \ref{j-coefficientsBigraded}, there exist elements $y_1, \ldots, y_d \in \m$ which  form an absolutely superficial sequence for both $\widetilde{T}$ and $\widetilde{T}/\widetilde{x}\widetilde{M}$.
By Propositions \ref{j-coefficientsBigraded} and \ref{j-multiplicityformula}, if we set $q=(y_1, \ldots, y_d)$, then for every $0\leq i \leq d-2$, we have

$$
(-1)^{i}j_{i}(I, \widetilde{M})=a_{(q,\, I,\, \widetilde{M})}^{(1,\,1)}(0,d-i)=a_{(q,\, I,\, \widetilde{M}/\widetilde{x}\widetilde{M})}^{(1,\,1)}(0,d-i)=(-1)^{i}j_{i}(I, \widetilde{M}/\widetilde{x}\widetilde{M}).
$$
\medskip
Furthermore for $1\leq i \leq d$, set $q_i=(y_1, \ldots, y_i)$, then
$$
[(y_1, \ldots, y_{i-1})\widetilde{T}:_{\widetilde{T}} y_{i}]\cap (y_1, \ldots, y_d)\widetilde{T}=(y_1, \ldots,y_{i-1})\widetilde{T}.
$$
$$
[(y_1, \ldots, y_{i-1})\widetilde{T}/\widetilde{x}\widetilde{M}:_{\widetilde{T}/\widetilde{x}\widetilde{M}} y_{i}]\cap (y_1, \ldots, y_d)\widetilde{T}/\widetilde{x}\widetilde{M}=(y_1, \ldots,y_{i-1})\widetilde{T}/\widetilde{x}\widetilde{M}.
$$
Since $\widetilde{T}$ and $\widetilde{T}/\widetilde{x}\widetilde{M}$ are finite modules over a polynomial ring over  the regular local ring $\widetilde{R}$, by \cite[3.6, 3.2]{T3}, there exists a dense open subset $U$ of $K^n$ such that if $\underline{\alpha}=(\alpha_1, \ldots, \alpha_n)\in U$, then
$$
[(y_1, \ldots, y_{i-1})T:_{T} y_{i}]\cap (y_1, \ldots, y_d)T
$$
$$
=([(y_1, \ldots, y_{i-1})\widetilde{T}:_{\widetilde{T}} y_{i}]\cap (y_1, \ldots, y_d)\widetilde{T})_{\underline{\alpha}}
$$
$$
=((y_1, \ldots,y_{i-1})\widetilde{T})_{\underline{\alpha}}=(y_1, \ldots,y_{i-1})T
$$
$$
[(y_1, \ldots, y_{i-1})T/xM:_{T/xM} y_{i}]\cap (y_1, \ldots, y_d)T/xM
$$
$$
=([(y_1, \ldots, y_{i-1})\widetilde{T}/\widetilde{x}\widetilde{M}:_{\widetilde{T}/\widetilde{x}\widetilde{M}} y_{i}]\cap (y_1, \ldots, y_d)\widetilde{T}/\widetilde{x}\widetilde{M})_{\underline{\alpha}}
$$
$$
= ((y_1, \ldots,y_{i-1})\widetilde{T}/\widetilde{x}\widetilde{M})_{\underline{\alpha}}=(y_1, \ldots,y_{i-1})T/xM.
$$
and the image $\overline{x}=\overline{\sum_{j=1}^n \alpha_{j}a_j}\in I/qI$ is filter-regular on $T^{\prime}={\rm gr}_{q}(T)$ with respect to the ideal $G^{\prime}_{(01)}G^{\prime}$.
Hence
$$
(-1)^{i}j_{i}(I, M)=a_{(q,\, I,\, M)}^{(1,\,1)}(0,d-i)=a_{(q,\, I,\, M/xM)}^{(1,\,1)}(0,d-i)=(-1)^{i}j_{i}(I, M/xM).
$$
for  every $0\leq i \leq d-2$.

\QED


  Lemma \ref{Singh} also allows us to keep track of the behaviors of $j_{d-1}$ and $j_d$ under a general hyperplane section as  it has been done for  ideals of definition on $M$.

It is well-known that the leading generalized Hilbert coefficient $j_0(I, M)$, which  is also called the $j$-multiplicity of $I$ on $M$,  is preserved under a general hyperplane section (see for instance \cite{NU}). The preservation of higher  generalized Hilbert coefficients was unknown due to the fact that it is hard to estimate the change of  the length of each homogenous component of $\Gamma_{\m}(T)$ under a general hyperplane section.  Our method gives another way to compute the length of homogenous components of $\Gamma_{\m}(T)$ by using a suitable bigraded module.

As we can see in Example \ref{ex2}  under a general hyperplane section, in general only the generalized coefficients  $j_0, \ldots, j_{d-2}$ are preserved. Indeed, let $S$ be the hypersurface ring obtained by the  polynomial ring $R=K[x,y,z,v,w]$ by modding out a general element $\xi$ in the ideal $I=I_2(B)$ as defined in Example \ref{ex2}. Let $L$ be the ideal $IS$. Notice that ${\rm dim} \, S=3$, the associated graded ring of $L$ is Cohen-Macaulay,  however $j_{2}(L)$ is not preserved if we go modulo a general element $\nu$ in $L$. Notice that $j_2(L/\nu L)=4>1= j_2(L)$, i.e. $j_{d-1}$ does not decrease if the initial form of the general element is regular on $G$. This is a general fact that is explained in the next Remark.

\begin{Remark}
Assume $R$ has infinite residue field and the grade of $G^+$ in $G$ is positive. For a general element $x$ in $I$ we have  $j_{d-1}(I, M/xM) \ge j_{d-1}(I,M)$. 
\end{Remark}

\noindent The  above remark follows directly from Proposition \ref{j-multiplicityformula} part (b), as the sum $$\sum_{\nu=2}^{t}\lambda(I^{\nu}M:_M x/
[(q^{s+1}I^{\nu}M+I^{\nu+1}M):_M x +I^{\nu-1}M])
$$ vanishes when the  initial form of  $x$ is regular on $G$.

\section{Applications of generalized Hilbert functions. }

\bigskip

In this section, we are going to discuss the generalized Hilbert function for ideals having minimal or almost minimal $j$-multiplicity on a finite module over a Noetherian local ring.

Let $(R,\m)$ be a Noetherian local ring, let $I$ be an $R$-ideal and let $M$ be a finite $R$-module of ${\rm dim}\,M=d$. Let $G$ and $T$ be defined as before. By adjoining a variable $z$ and  passing to the local ring $R[z]_{\m R[z]}$, we may assume that the residue field of $R$ is infinite. Recall for general elements $x_1, \ldots, x_d$ in $I$, let $\overline{M}=M/(x_1,\ldots,x_{d-1})M:_M I^{\infty}$. Then $\overline{M}$ is either the zero module or a 1-dimensional Cohen-Macaulay module and $I$ is an ideal of definition on $\overline{M}$.
Then one has the Hilbert-Samuel function
of $I$ on $\overline{M}$
$$
H_{I,\,\overline{M}}(t)=\lambda_{\overline{R}/\overline{I}}(\overline{M}/I^{t+1}\overline{M})
$$
and the Hilbert-Samuel series of $I$ on $\overline{M}$
$$
h_{I,\,\overline{M}}(z)=\sum_{t\geq 0}H_{I,\,\overline{M}}(t)z^t.
$$

The following lemma says that the Hilbert-Samuel function $H_{I, \,\overline{M}}$  does not depend on the general elements $x_1, \ldots, x_d$ in $I$. The  proof follows from  the ingredient proved in  \cite[2.3]{PX}, where the independence of  $\lambda(I\overline{M}/I^2\overline{M})$ is established.

\begin{Lemma}
Let $R$ be a Noetherian local ring with infinite residue field, let $I$ be any $R$-ideal, and let $M$ be a finitely generated $R$-module of dimension $d$.   Then for general elements $x_1, \ldots, x_d$ in $I$, $\lambda(I^t\overline{M}/I^{t+1}\overline{M})$ does not depend on $x_1, \ldots, x_d$ for every $t\geq 0$.
\end{Lemma}

\demo
If $\ell(I, M)<d$ then $\overline{M}=0$. Thus we may assume that $\ell(I, M)=d$. By \cite[2.5]{PX}, $j(I, M)=j(I, I^tM)$ for every $t\geq 0$. Therefore the case where $t\geq 1$ follows from \cite[2.3]{PX}. The case where $t=0$ also follows by a similar argument as in the proof of \cite[2.3]{PX}.
\QED

\vspace{0.3in}

Recall if $I$ has minimal $j$-multiplicity on $M$, by \cite[Theorem 2.9]{RV},
$$
h_{I, \,\overline{M}}(t)=\frac{h_0+h_1t}{(1-t)^2}
$$
where
$$h_0=\lambda(\overline{M}/I\overline{M})=\lambda(M/[(x_1,\ldots,x_{d-1})M:_M I^{\infty}+IM])$$
 and
 $$h_0+h_1=\lambda(I\overline{M}/I^2\overline{M})=\lambda(I M/[(x_1,\ldots,x_{d-1})M:_{IM} I^{\infty}+I^2M]).$$

 If $I$ has almost minimal $j$-multiplicity on $M$, by \cite[Corollary 4.4]{RV},
$$
h_{I,\, \overline{M}}(t)=\frac{h_0+h_1t+t^s}{(1-t)^2}
$$
for some $s\geq 2$, where  $h_0$ and $h_1$ are as above.

Furthermore if $I$ is an ideal of definition on $M$, the Hilbert function is preserved after modding out
elements which are regular on the associated graded module $T$. Hence the shape of the Hilbert-Samuel series of $M$ is the same as the Hilbert-Samuel series of $\overline{M}$. However this is no longer true for arbitrary ideals as the following two examples show.

\begin{Example}
Let $R=k[x,y,z,v,w]$ be a polynomial ring over an infinite field $k$ and $I=I_2(A)$, where
$$
A=\left(\begin{matrix}
x & y & z & v \\
y & z & v & w
\end{matrix}\right).
$$
Observe that $I$ is a perfect ideal of grade 3,  $I$ is a complete intersection on the punctured spectrum with analytic spread $\ell(I)=5$. Hence $I$ satisfies $G_5$ and $AN_5^-$. By  Macaulay 2, $\lambda(I^2/JI)=0$, where $J$ is a general minimal reduction of $I$. Hence $I$ has minimal $j$-multiplicity. By \cite[4.9]{PX}, the associated graded ring $G={\rm gr}_I(R)$ is Cohen-Macaulay.

However by Macaulay 2, the generalized Hilbert-Samuel series of $I$ is:
$$h_I(z)=\frac{z+z^2+z^3+z^4}{(1-z)^6}$$
which does not have the expected shape.

After modding out a general element $\xi_1$, two general elements $\xi_1, \xi_2$, three general elements $\xi_1, \xi_2, \xi_3$, and  four general elements $\xi_1, \xi_2, \xi_3, \xi_4$ in $I$, the generalized Hilbert-Samuel series are $\frac{z+z^2+z^3+z^4}{(1-z)^5}$, 
 $\frac{z+z^2+z^3+z^4}{(1-z)^4}$, $\frac{z+z^2+z^3+z^4}{(1-z)^3}$ and $\frac{4z}{(1-z)^2}$ respectively.
 The Hilbert-Samuel series modulo the residual intersection $(\xi_1, \xi_2, \xi_3, \xi_4) :I^{\infty}$ is $h_{\overline{I}}(z)=\frac{3+z}{(1-z)^2}$.
\end{Example}

In the previous example the generalized Hilbert series was preserved until the height of the ideal was reduced to zero. One could imagine that this is always true if $G$ is Cohen-Macaulay. However the next example shows that the Cohen-Macaulayness of $G$ does not suffice to preserve the generalized Hilbert-Samuel series:

\begin{Example}\label{ex2}
Let $R=k[x,y,z,v]$ be a polynomial ring over an infinite field $k$ and $I=I_2(B)$, where
$$
B=\left(\begin{matrix}
x & y & z & v \\
v & x & y & z
\end{matrix}\right).
$$
Observe that $I$ is a perfect ideal of grade 3, $I$ is a generically a complete intersection with analytic spread $\ell(I)=4$. Hence $I$ satisfies $G_4$ and $AN_4^-$. 
Let $J$ be a general minimal reduction of $I$, using Macaulay 2, one computes $\lambda(I^2/JI)=1$. 
Hence $I$ has almost minimal $j$-multiplicity. By \cite[4.9]{PX}, the associated graded ring $G={\rm gr}_I(R)$ has depth at least $3$. Indeed, $G$ is Cohen-Macaulay.

However by Macaulay 2, the generalized Hilbert-Samuel  series of $I$ is:
$$h_I(z)=\frac{4z+z^2+6z^3-3z^4}{(1-z)^5}$$
which does not have the expected shape.

After modding out a general element $\xi_1$, two general elements $\xi_1, \xi_2$, and three general elements $\xi_1, \xi_2, \xi_3$ in $I$ , the generalized Hilbert-Samuel series are $\frac{4z+z^2+6z^3-3z^4}{(1-z)^4}$, $\frac{4z+4z^2}{(1-z)^3}$ and   $\frac{7z+z^2}{(1-z)^2}$ respectively. The Hilbert-Samuel series modulo the residual intersection $(\xi_1, \xi_2, \xi_3) :I^{\infty}$ is $h_{\overline{I}}(z)=\frac{1+6z+z^2}{(1-z)^2}$.
\end{Example}

Finally we provide a sufficient  condition to ensure that the generalized Hilbert-Samuel series of ideals having minimal or almost minimal multiplicity  has the desired shape.

\begin{Theorem}\label{Minimal}
Let $M$ be a  Cohen-Macaulay module of dimension $d$ over a Noetherian local ring $(R, \m)$ and let $I$ be an $R$-ideal with $\ell(I, M)=d$. Assume  ${\rm depth}\,(M/IM)\geq 1$ and  $I$ satisfies $G_{d}$ and  $AN^-_{d-2}$ on $M$. Let $q\subseteq \m$ be an ideal of definition on $M/IM$ that is generated by a $d$-sequence on the associated graded module $T$. Let  $x_1, \ldots, x_d$ be general elements in $I$ and set  $J_i=(x_1, \ldots, x_i)$ for $1\leq i\leq d$.
\begin{itemize}
\item[(a)] Assume $I$ has minimal $j$-multiplicity and for $1\leq i \leq d-1$ and $s>>0$,
$$
J_iIM\cap (q^{s+1}I IM+I^{2}IM)=J_i(q^{s+1}IM+I^2M),
$$
 then  $h_{I, \,M}(z)=\frac{\varepsilon_0 z+ \varepsilon_1 z^2}{(1-z)^{d+1}}$,
 where $\varepsilon_0=\lambda(\Gamma_{\m}(IM/I^2M))$ and $\varepsilon_1=h_0+h_1-\varepsilon_0$ with 
 $h_0+h_1=\lambda(I\overline{M}/I^2\overline{M})=\lambda(IM/(x_1, \ldots, x_{d-1})M+I^2M)$  defined above.

\item[(b)] Assume $I$ has almost minimal $j$-multiplicity and for $1\leq t\leq r(I, M),$ $1\leq i \leq d-1$ and 
$s>>0$,
$$
J_iI^{t-1}IM\cap (q^{s+1}I^t IM+I^{t+1}IM)=J_i(q^{s+1}I^{t-1}IM+I^t IM),
$$ 
then   $h_{I, \,M}(z)=\frac{\varepsilon_0 z+\varepsilon_1 z^2+z^s}{(1-z)^{d+1}}$ for some integer $s\geq 3$,
where $\varepsilon_0$ and $\varepsilon_1$ are the same as in part (a).
\end{itemize}
\end{Theorem}

\demo By \cite{PX}, $T$ is Cohen-Macaulay (respectively, almost Cohen-Macaulay) if $I$ has minimal (respectively, 
almost minimal) $j$-multiplicity. Furthermore the initial forms $x_1^*, \ldots, x_{d-1}^*$ in $G$ form a regular sequence on $T_+$. Hence $I^{t}IM:_{IM} x= I^{t-1}IM$ for every $t\geq 1$. Appling Lemma \ref{Singh} to the module $IM$, one has for every $s, t \geq 0$,
$$\vspace{0.07in}
\begin{array}{l}
H_{(q,\, I,\, IM)}^{(1,\,0)}(s,t)=H_{(q,\, I,\, IM/x_1IM)}^{(1,\,1)}(s,t)\vspace{0.1in}\\
-\sum_{\nu=1}^t \lambda([(q^{s+1}I^\nu IM+I^{\nu+1}IM):_{I^{\nu-1}IM}x]/[q^{s+1}I^{\nu-1}IM+I^\nu IM]).
\end{array}\vspace{0.07in}
$$

By \cite{PX}, $r(I, M)=1$ in part (a). 
Applying a similar argument as in \cite{CGPU}, one can show that for $s>>0$, the images of $x_1, \ldots, x_{d-1}$ in $T_+/q^{s+1}T_+$ form a regular sequence on $T_+/q^{s+1}T_+$. Hence
$(q^{s+1}I^t IM+I^{t+1}IM):_{I^{t-1}IM}x=q^{s+1}I^{t-1}IM+I^t IM$ for every $t\geq 1$.
Therefore 
$H_{(q,\, I,\, IM)}^{(1,\,0)}(s,t)=H_{(q,\, I,\, IM/x_1IM)}^{(1,\,1)}(s,t)$ for every $s, t \geq 0$.
This shows that the h-polynomial of the generalized Hilbert-Samuel series of $I$ on $IM$ are preserved after modulo $x_1, \ldots, x_{d-1}$.
Since  $\Gamma_0(M/IM)=0$, one has $H_{(q,\, I,\, M)}^{(1,\,0)}(s,t)= H_{(q,\, I,\, IM)}^{(1,\,0)}(s,t-1)$ and $h_{I, \,M}(z)=h_{I, \,IM}(z)z$.     After lifting back the Hilbert-Samuel series in the 1-dimensional case, we obtain the desired result.
\QED



\end{document}